\documentclass[12pt,leqno]{article}

\usepackage{amssymb}
\usepackage{amsmath}
\usepackage{color}
\usepackage{url}

\newcommand{\V}{\mbox{\sf V}}
\newcommand{\F}{\mbox{$\mathfrak F$}}
\newcommand{\A}{\mbox{$\mathfrak A$}}
\newcommand{\B}{\mbox{$\mathfrak B$}}
\newcommand{\al}{\alpha}
\newcommand{\be}{\beta}
\newcommand{\ga}{\gamma}
\newtheorem{thm}{Theorem}

\begin{document}

\title{Free algebras of discriminator varieties generated by finite algebras are atomic}
\author{Andr\'eka, H.\ and N\'emeti, I.}%
\date{February 2016}
\maketitle

\begin{abstract} We prove that all definable pre-orders are atomic, in a
finitely generated free algebra of a discriminator variety of finite
similarity type which is generated by its finite members.
\end{abstract}

A pre-order $\le$ is a reflexive and transitive relation, $a<b$
abbreviates that $a\le b$ and $b\not\le a$. A pre-order $\le$ is
called \emph{atomic} if in each interval the smallest element has a
cover, formally if $\forall a<b\,\exists c\,[a<c\le
b\land\neg\exists x(a<x<c)]$. The name is justified by the fact that
a Boolean algebra is atomic iff its natural order is atomic in the
above sense. A binary relation $\le$ in an algebra $\F$ is called
\emph{definable} if it is equationally definable, i.e., if there are
terms $\tau,\sigma$ in the language of $\F$ such that
$\F\models\forall xy[x\le y\leftrightarrow \tau=\sigma]$.
A variety $\V$ is a \emph{discriminator variety} if there is a term
$\sigma$ in its language that in each subdirectly irreducible member
of $\V$ is the so-called \emph{switching term}, i.e.,
$\sigma(x,y,u,v)=u$ if $x=y$ and $\sigma(x,y,u,v)=v$ if $x\ne y$.
For this definition of a discriminator variety, and for its basic
properties that we will use in the proof, we refer the reader to
\cite{UA}.
Theorem~\ref{atomic} below is a corollary of \cite[Theorem
4.1(i)]{AJN91} which states that in an ordered discriminator variety
of finite similarity type, each finitely generated residually finite
algebra is atomic. (The quoted theorem concerns only orders and not
pre-orders, but its proof equally applies to pre-orders.) Here, we
give a different, direct proof that might be easier to generalize in
certain directions.

\begin{thm}\label{atomic}
Let $\V$ be a discriminator variety of finite similarity type that
is generated by its finite members as a variety. Any definable
pre-order in a finitely generated free algebra of $\V$ is atomic.
\end{thm}

We note that none of the conditions of the theorem can be omitted
without affecting its truth, see the discussion at the end of the
paper. Theorem~\ref{atomic} has applications in logic, see
\cite{Kha16}. Both its statement and proof are a generalization of
\cite[Thm.2.5.7]{HMT}. %Another generalization that follows the
%original one more closely can be found in \cite{Sayed}.
The key ingredients of our proof are that all pre-orders on a finite
set are atomic and the discriminator term gives us expressive
power.\bigskip

\noindent{\bf Proof of Theorem~\ref{atomic}.}
Let $\F$ be a $\V$-free algebra freely generated by a finite set
$X$. Assume that the definable binary relation $\le$ in $\F$ is a
pre-order. Let $\al < \be$ in $\F$, i.e., $\al\le\be$ and $\be\not
\le \al$. We set to finding a cover $\ga\le \be$ of $\al$.

Since $\al < \be$ in $\F$ and $\V$ is generated by finite algebras,
there are a finite subdirectly irreducible $\A\in\V$ and a
homomorphism $h:\F\to \A$ such that $h(\al)<h(\be)$ (in $\A$). We
may assume that $\A$ is generated by $h(X)=\{ h(x) : x\in X\}$ since
a subalgebra of a subdirectly irreducible algebra in a discriminator
variety is also subdirectly irreducible. Since $X$ generates $\F$,
each element of $F$ is of the form $\tau^{\F}(\bar{x})$ for at least
one term $\tau$ (where $\bar{x}$ is a sequence of elements from
$X$). Since $\A$ is also generated by $h(X)$, the homomorphism $h$
is surjective. For each $a\in A$ let $\rho(a)$ denote a term such
that $h(\rho(a))=a$. ($\rho$ stands for ``representative".)
Consider the formulas of the following three forms
$$x=\rho(h(x)),\quad f(\rho(a_1),\dots \rho(a_m))=\rho(f(a_1,\dots,a_m)),
\quad \rho(a)\ne\rho(a')$$ where $x\in X$, $f$ is an $m$-place
operation in $\A$, $a_1,\dots,a_m\in A$, and $a,a'\in A$ are
distinct. These conditions specify a finite set of formulas because
$X$, the similarity type of $\A$, and $A$ are all finite.  Thus
their conjunction $\eta$ is a universal formula with $x\in X$
considered to be free variables.
% (note that all the terms in $F$ are written up from elements of $X$).
Let $\A\models\eta[h]$ denote that $\eta$ is true in $\A$ under the
evaluation $h$ of variables, and the same for $\B, k$. The formula
$\eta$ describes $\A$ up to isomorphism, but we will use only
$\A\models\eta[h]$ and the following property of $\eta$,
\begin{description}
\item[(1)]\label{eq}
Let $k:\F\to\B$ be such that $\B\models\eta[k]$. Then\\ $i=\{
\langle h(\tau),k(\tau)\rangle : \tau\in F\}$ is an embedding of
$\A$ into $\B$.
\end{description}
Indeed, recall that $X$ generates $\F$ and so each element of $F$
can be obtained from $X$ by application of a term $\tau$. First one
shows by induction on $\tau$ that $k(\tau)=k(\rho(h(\tau)))$ for all
$\tau$. Indeed, we have $k(x)=k(\rho(h(x)))$ by $\eta$, and we have
$k(f(\tau_1,\dots,\tau_m))=f(k(\tau_1),\dots,k(\tau_m))=f(k(\rho(h\tau_1)),\dots,k(\rho(h\tau_m)))=$\\
= $k(f(\rho(h\tau_1),\dots,\rho(h\tau_m)))=
k(\rho(f(h\tau_1,\dots,h\tau_n)))= k(\rho(hf(\tau_1,\dots,\tau_m))$
by $k$ being a homomorphism, the induction hypothesis, $k$ being a
homomorphism, $\eta$, and $h$ being a homomorphism, respectively.
Thus $i(a)=k\rho(a)$ for all $a\in A$. Hence $i:\A\to\B$ is a
homomorphism by $k$ being a homomorphism and the second kind of
terms in the definition of $\eta$. It is one-to-one by the last term
in the definition of $\eta$ because $a=h(\rho(a))$ for all $a\in A$.

Since we are in a discriminator variety, each universal formula is
equivalent to an equation in the subdirectly irreducible algebras.
Let $\delta,\varepsilon$ be terms such that the formula
$$ \forall\bar{x}[\eta(\bar{x})\leftrightarrow\delta=\varepsilon]$$
is valid in the subdirectly irreducible members of $\V$. Recall that
$h(\al)<h(\be)$ in $\A$. Since $\A$ is finite, there is a cover
$c\le h(\be)$ of $h(\al)$ in $\A$. Now, we are ready to write up the
desired cover term $\gamma$ of $\al$. Let $\gamma$ denote
$$\sigma(\delta,\varepsilon,\rho(c),\al)$$
where $\sigma$ is the switching term of the discriminator variety
$\V$.
We are going to show that $\gamma\le\be$ is a cover of $\al$ in
$\F$.

We begin with showing that $\al \le \gamma\le\be$ in $\F$. In order
to show this, we have to show $\al\le\gamma\le\be$ in all
subdirectly irreducible elements $\B$ of $\V$. Notice that
$\al\le\rho(c)$ does not necessarily hold in $\F$. Let thus
$k:\F\to\B\in\V$ be arbitrary, with $\B$ subdirectly irreducible.
Assume first that $k(\delta)=k(\varepsilon)$. Since $k$ is a
homomorphism,
$k(\gamma)=k(\sigma(\delta,\varepsilon,\rho(c),\al))=\sigma(k(\delta),k(\varepsilon),k(\rho(c)),k(\al))$.
Since $\sigma$ is a switching term in $\B$ and
$k(\delta)=k(\varepsilon)$, we get $k(\gamma)=k(\rho(c))$ and
$\B\models\eta[k]$.
By (1) and $h(\al)\le h(\rho(c))\le h(\be)$ then $k(\al)\le
k(\rho(c))\le k(\be)$ and we are done. Assume that $k(\delta)\ne
k(\varepsilon)$. Then $k(\gamma)=k(\al)$ by the definition of
$\gamma$, so again we are done since $\al\le\be$ in $\F$.

Next we show $\gamma\not\le\al$ in $\F$. By the construction of
$\eta$ we have that $\A\models\eta[h]$. Since $\A$ is subdirectly
irreducible then $\A\models(\delta=\varepsilon)[h]$, i.e.,
$h(\delta)=h(\varepsilon)$ by the choice of $\delta,\varepsilon$.
Thus $h(\gamma)=h(\rho(c))$ by the definition of $\gamma$, but
$h(\rho(c))=c$ by the choice of $\rho(c)$. This shows that
$h(\gamma)\not\le h(\al)$ since $h(\al)<c$ in $\A$, and this implies
$\gamma\not\le\al$ in $\F$.

Finally, we show that $\gamma$ is a cover of $\al$ in $\F$. We will
show that for all $\tau\in F$ such that $\al \le \tau \le \gamma$ in
$\F$ we have either $\tau\le\al$ or $\gamma\le\tau$ in $\F$. So,
assume $\al \le \tau \le \gamma$.
Then $h(\al)\le h(\tau)\le h(\gamma)$. Since $h(\gamma)=c$ is a
cover of $h(\al)$ in $\A$, either $h(\tau)\le h(\al)$, or
$h(\gamma)\le h(\tau)$. Let $k:\F\to\B$ be arbitrary with $\B\in\V$
subdirectly irreducible.
Assume $h(\tau)\le h(\al)$. If $k(\delta)\ne k(\varepsilon)$ then
$k(\gamma)=k(\al)$, so $k(\tau)\le k(\al)$ by $\tau\le\gamma$. If
$k(\delta)=k(\varepsilon)$, then $k(\tau)\le k(\al)$ by
$\B\models\eta[k]$, (1) and $h(\tau)\le h(\al)$. So, in either case
we have $k(\tau)\le k(\al)$, which implies that $\tau\le\al$ in
$\F$.
Assume $h(\gamma)\le h(\tau)$. If $k(\delta)\ne k(\varepsilon)$ then
$k(\gamma)=k(\al)$, so $k(\gamma)\le k(\tau)$ by $\al\le\tau$. If
$k(\delta)=k(\varepsilon)$ then $\B\models\eta[k]$, so by (1) we
have $k(\gamma)\le k(\tau)$ by $h(\gamma)\le h(\tau)$. Since in
either case $k(\gamma)\le k(\tau)$, we have $\gamma\le\tau$ in $\F$,
and we are done.
%
%Let now $k:\F\to\B$ be arbitrary with $\B\in\V$ subdirectly
%irreducible. If $k(\delta)\ne k(\varepsilon)$ then $k(\tau)\le
%k(\al)$ by $k(\gamma)=k(\al)$ and $\tau\le\gamma$. Assume
%$k(\delta)=k(\varepsilon)$. Then $\B\models\eta[k]$, so by (1) we
%have $k(\tau)=ih(\tau)$ and $k(\al)=ih(\al)$.  Assume that
%$h(\tau)\le h(\al)$. Then $k(\tau)\le k(\al)$ for all $k:\F\to\B$
%which means that $\tau\le \al$ in $\F$. Assume that $c\le h(\tau)$.
%Then $k(\tau)=k(\al)$ when $k(\delta)\ne k(\varepsilon)$ and
%$k(\tau)=k(\rho(c))$ otherwise. This means that
%$k(\tau)=\sigma(k(\delta),k(\varepsilon),k(\al),k(\rho(c)))$=$k(\sigma(\delta,\varepsilon,\al,\rho(c))$
%=$k(\gamma)$ for all $k$, so $\tau=\gamma$ in $\F$, and we are done.
\hfill {\sf QED}\bigskip

We have the following corollary.

\begin{thm}\label{BA}
Let $\V$ be a discriminator variety that is generated by its finite
algebras as a variety. Assume that $\V$ has a Boolean algebra reduct
and its similarity type is finite. Then the finitely generated
$\V$-free algebras are atomic.
\end{thm}

Most varieties arising from logic have Boolean algebra reducts, and
atomicity of their free algebras corresponds to weak G\"odel's
incompleteness property holding for the logic, see \cite{NPrep85} or
\cite{Gyenis, Kha16}.

None of the conditions of Theorem~\ref{atomic} can be omitted
without affecting its truth. The condition of $\V$ being of finite
similarity type is necessary because the free algebra of the class
of Boolean algebras with infinitely many constants (of which we do
not state any equations) is atomless. The condition of the free
algebra generated by finitely many elements is necessary because the
infinitely generated free Boolean algebra is atomless. The condition
that $\V$ is discriminator is necessary because the variety of
2-dimensional cylindric-relativized set algebras is generated by its
finite members but its finitely generated free algebras are not
atomic, see \cite{Kha16}. The case is similar for other varieties of
relativized algebras in algebraic logic, see \cite{Kha15a, Kha15b}.
The condition that $\V$ be generated by its finite members is
necessary, because varieties of un-relativized algebras in algebraic
logic usually do not have atomic free algebras but they are
discriminator in the finite-dimensional case. This is the case for
the varieties of abstract and representable relation algebras, the
abstract and representable finite-dimensional cylindric algebras,
diagonal-free 3-dimensional cylindric algebras. See, e.g.,
\cite{Gyenis, HH02, NPrep85, TG} or \cite[4.3.32]{HMT}.

\bigskip\bigskip\bigskip

\noindent Alfr\'ed R\'enyi Institute of Mathematics, Hungarian
Academy of Sciences\\
Budapest, Re\'altanoda st.\ 13-15, H-1053 Hungary\\
andreka.hajnal@renyi.mta.hu, nemeti.istvan@renyi.mta.hu
\end{document}